\documentclass[11pt]{amsart}

\usepackage{amsthm,amsmath,amssymb,amscd}
\usepackage{amssymb}
\usepackage{graphics}
\newtheorem{theorem}{Theorem}[section]
\newtheorem{corollary}[theorem]{Corollary}

\newtheorem{lemma}[theorem]{Lemma}
\newtheorem{proposition}[theorem]{Proposition}
\newtheorem{definition}{Definition}[section]

\numberwithin{equation}{section}

\begin{document}

\title [Generalized Lagrangian mean curvature flows]
 {Generalized Lagrangian mean curvature flows in almost Calabi-Yau manifolds}
\author{Jun Sun, Liuqing Yang}

\address{CIRM, Fondazione Bruno Kessler, Via Sommarive, 14 - Povo, I-38123 \\ Trento (TN), Italy.}
\email{sunjun@fbk.eu}

\address{Academy of Mathematics and Systems Sciences, Chinese Academy of Sciences, Beijing 100190, P. R. of China.}
\email{yangliuqing@amss.ac.cn}


\date{}

\maketitle
\begin{abstract}
In this paper, we study the generalized Lagrangian mean curvature
flow in almost Einstein manifold proposed by T. Behrndt. We show
that the singularity of this flow is characterized by the second
fundamental form. We also show that the rescaled flow at a
singularity converges to a finite union of Special Lagrangian cones
for generalized Lagrangian mean curvature flow with zero-Maslov
class in almost Calabi-Yau manifold. As a corollary, there is no
finite time Type-I singularity for such a flow.
\end{abstract}

\vspace{.2in}


\section{Introduction}

\allowdisplaybreaks

\vspace{.1in}

\noindent Suppose $(M,J,\bar{\omega},\bar{g})$ is a smooth K\"ahler
manifold with complex dimension $n$, complex structure $J$, K\"ahler
metric $\bar{g}$ and K\"ahler form $\bar{\omega}$. The K\"ahler form
and the K\"ahler metric are related by
\begin{equation*}
    \bar{\omega}(X,Y)=\bar{g}(JX,Y),
\end{equation*}
for $X,Y \in \Gamma(TM)$. Moreover, suppose that $\overline{Ric}$ is
the Ricci tensor of $\bar{g}$. Then the Ricci form $\bar{\rho}$ is
defined by
\begin{equation*}
    \bar{\rho}(X,Y)=\overline{Ric}(JX,Y).
\end{equation*}
Let $L$ be a compact manifold of real dimension $n$ and $F_0:
L\longrightarrow M$ an immersion of $L$ into $M$. The induced metric
on $L$ is $g=F_0^*\bar{g}$ and set $\omega=F_0^*\bar{\omega}$. It is
known by definition that $F_0$ is a \emph{Lagrangian immersion} if
$\omega=0$.

In 1996, Strominger, Yau and Zaslow (\cite{SYZ}) found that mirror
symmetry is related to special Lagrangian submanifold (which is automatically minimal) in Calabi-Yau
manifold. One natural approach to obtaining minimal submanifold is
to evolve a submanifold along the negative gradient flow of the area
functional, i.e., the mean curvature flow. Fortunately, when the
ambient manifold $M$ is K\"ahler-Einstein, Smoczyk (\cite{Sm1})
proved that if the initial surface $L_0$ is Lagrangian, then along
the mean curvature flow, it remains Lagrangian for each time. Since
then, Lagrangian mean curvature flow received a lot of attention and
there are many results on it. (c.f. \cite{CL2}, \cite{N}, \cite{TY},
\cite{Wa}, etc.) All of them concern Lagrangian mean curvature flow
in K\"ahler-Einstein manifold, while most of them focus on
Calabi-Yau ambient manifold.

Recently, generalized Lagrangian mean curvature flow attracts more
attention (\cite{B}, \cite{SW}). This flow was first studied by T.
Behrndt (\cite{B}). Instead of considering mean curvature flow in a
K\"ahler-Einstein manifold, he considered the case when the ambient
manifold is almost Einstein. Let us first recall the definition of
an almost Einstein manifold in \cite{B}.
\begin{definition}
An $n$-dimensional K\"ahler manifold $(M,J,\bar{\omega},\bar{g})$ is
called almost Einstein if
\begin{equation*}
    \bar{\rho}=\lambda\bar{\omega}+ndd^c\psi
\end{equation*}
for some constant $\lambda\in \textbf{R}$ and some smooth function
$\psi$ on $M$.
\end{definition}
Suppose the  K\"ahler manifold $(M,J,\bar{\omega},\bar{g})$ is
almost Einstein. Given an immersion $F_{0}: L \longrightarrow M$ of
a manifold $L$ into $M$, T. Behrndt (\cite{B}) proposed the
generalized mean curvature flow,
\begin{equation}\label{flow}
   \begin{cases}
     \frac{\partial}{\partial t} F(x,t) =\textbf{K}(x,t), \ \ (x,t)\in L\times (0,T) \\
      F(x,0)=F_{0}(x), \ x\in L.
   \end{cases}
\end{equation}
Here \begin{equation*}
    \textbf{K}=\textbf{H}-n\pi_{\nu L}(\bar{\nabla}\psi)
\end{equation*}
is a normal vector field along $L$ which is called the generalized
mean curvature vector field of $L$. As \textbf{K} is a differential
operator differing from \textbf{H} just by lower order terms, it is
easy to see that (\ref{flow}) has a unique solution on a short time
interval (\cite{B}).

Arguing in a similar way as Smoczyk did for K\"ahler-Einstein case
(\cite{Sm1}), Behrndt (\cite{B}) proved that if $L_0=F_0(L)$ is
Lagrangian in the almost Einstein manifold $M$, then along the
generalized mean curvature flow (\ref{flow}), it remains Lagrangian
for each time. Therefore, it is reasonable to call such a flow
generalized Lagrangian mean curvature flow.

As a special case, Behrndt (\cite{B}) also considered the
generalized Lagrangian mean curvature flow in an almost Calabi-Yau
manifold. Let us recall the definition of an almost Calabi-Yau
manifold in \cite{J}.
\begin{definition}
An $n$-dimensional almost Calabi-Yau manifold
$(M,J,\bar{\omega},\bar{g},\Omega)$ is an $n$-dimensional K\"ahler
manifold $(M,J,\bar{\omega},\bar{g})$ together with a non-vanishing
holomorphic volume form $\Omega$.
\end{definition}

It can be seen that (\cite{B}), there exists a smooth function
$\psi$ on an almost Calabi-Yau manifold $M$ such that the Ricci form
of $(M,\bar{g})$ is given by
\begin{equation*}
    \bar{\rho}=ndd^c\psi.
\end{equation*}
In particular, this implies that an almost Calabi-Yau manifold is
almost Einstein.

Similar to the Calabi-Yau case, we can define the Lagrangian angle
$\theta:L \to {\textbf{S}}^1$ for a Lagrangian submanifold in an
almost Calabi-Yau manifold, which satisfies (\cite{HarLaw})
\begin{equation*}
    F^*\Omega=e^{i\theta+n\psi}d\mu_g
\end{equation*}
for $F:L\to M$ a Lagrangian immersion. Note that $\theta$ is a
multi-valued function on $L$, which is well-defined up to an
additive constant $2k\pi$, $k\in\textbf{Z}$. Behrndt (\cite{B})
proved that on a Lagrangian submanifold of an almost Calabi-Yau
manifold, we have
\begin{equation}\label{e1.2}
    \textbf{K}=J\nabla \theta.
\end{equation}
Furthermore, along the generalized Lagrangian mean curvature flow,
the Lagrangian angle satisfies (Proposition 5 of \cite{B})
\begin{equation}\label{e1.3}
      \frac{\partial}{\partial t}\theta=\Delta\theta+nd\psi(\nabla\theta).
\end{equation}
We call a Lagrangian submanifold \emph{almost calibrated} if
$\cos\theta>0$. When the Lagrangian angle $\theta$ is a single
valued function, the Lagrangian $L$ is called \emph{zero-Maslov}. By
(\ref{e1.3}) one can easily show that zero-maslov condition is
preserved under the generalized Lagrangian mean curvature flow
(\ref{flow}). It is obvious that almost calibrated Lagrangian must
be zero-Maslov class. We call a Lagrangian submanifold $L$
\emph{Special Lagrangian} if $\theta\equiv \theta_0$ is a constant
fucntion on $L$ (see Definition 5 and Proposition 3 of \cite{B}). In
this case, $L$ is calibrated with respect to
$Re(e^{-i\theta_0}\Omega)$ for the metric $\tilde{g}\equiv
e^{2\psi}\bar{g}$.

Recall that $L$ is a Lagrangian submanifold if $\bar{\omega}|_{L}\equiv0$. Likewise, as in \cite{N}, we define an integral
$n$-varifold $L_1$ and an integral $n$-current $L_2$ to be
\emph{Lagrangian} if
\begin{equation*}
    \int_{L_1} \phi |\omega\wedge\eta|d\mu=0 \ \ for\ all\ n-2 \ form \ \eta \ and\ all\ smooth\ \phi\in C^{\infty}(M)
\end{equation*}
and
\begin{equation*}
    \int_{L_2} \phi \omega\wedge\eta=0 \ \ for\ all\ n-2 \ form \ \eta \ and\ all\ smooth\ \phi\in C^{\infty}(M)
\end{equation*}
respectively. The concept of being Special Lagrangian can be easily
extended to the case when $L$ is an integral current.

\vspace{.1in}

It is known that, the mean curvature flow will blow up as the
maximal norm of the second fundamental form blows up. According to
the blow up rate, Huisken (\cite{Hu2}) divided the singularities of
mean curvature flow into two types: Type-I and Type-II. Generally,
singularity of mean curvature flow is unavoidable. Smoczyk (Theorem
2.3.5 of \cite{Sm2}) first proved that there is no compact Type-I
singularities with zero-maslov class. Later, Chen-Li (\cite{CL2})
and Wang (\cite{Wa}) independently proved that there is no Type-I
singularity for almost calibrated Lagrangian mean curvature flow in
Calabi-Yau manifold. Recently, Neves (\cite{N}) proved that there is
no finite time Type-I singularity for Lagrangian mean curvature flow
with zero-Maslov class. On the contrary, in 2007,
Groh-Schwarz-Smoczyk-Zehmisch (\cite{GSSZ}) constructed examples of
monotone, equivariant Lagrangian mean curvature flow which can
develop Type-I singularity.

Motivated by the previous work on Lagrangian mean curvature flow in
Calabi-Yau manifold, in this paper, we will study the generalized
Lagrangian mean curvature flow in almost Einstein manifold and
almost Calabi-Yau manifold. By computing the evolution equation of
the second fundamental form, we can see that (Theorem 2.4) the
blowing up of the generalized Lagrangian mean curvature flow
(\ref{flow}) in an almost Einstein manifold is also characterized by
the maximal norm of the second fundamental form.

\vspace{.1in}

Assume now that $M$ is an almost Calabi-Yau manifold and the solution to generalized Lagrangian mean
curvature flow develops a singularity at the point $(X_0,T)$ in
space-time. We consider the rescaled flow
\begin{equation}\label{recaled}
    F_{\lambda}(x,t)=\lambda(F(x,T+\lambda^{-2}s)-X_0), \ \ \ for \ -\lambda^2 T<s<0.
\end{equation}
We denote the scaled submanifold by
$(L^{\lambda}_{s},d\mu^{\lambda}_{s})$. Given any $\{\lambda_i\}$
going to infinity, we denote $L^{\lambda_i}_{s}$ by $L^{i}_{s}$ . By
establishing a monotonicity formula and following the argument of
Neves (\cite{N}), we can prove that

\vspace{.1in}

\noindent {\bf Main Theorem} {\it Let
$(M,J,\bar{\omega},\bar{g},\Omega)$ be a compact almost Calabi-Yau
manifold of complex dimension $n$. Assume that the initial surface
$L_0$ is Lagrangian and zero-Maslov class. Then for any sequence of
rescaled flows $(L_i^s)_{s<0}$ at a singularity there exist a finite
set $\{\bar{\theta}_1, \cdots, \bar{\theta}_N\}$ and integral Special
Lagrangian cones
$$L_1, ..., L_N$$
such that, after passing to a subsequence, we have for every smooth
function $\phi$ compactly supported, every $f$ in $C^2(\textbf{R})$,
and every $s<0$
\begin{equation*}
    \lim_{i\to\infty}\int_{L^{i}_{s}}f(\theta_{i,s})\phi d\mu^{i}_{s}=\sum_{j=1}^{N}m_jf(\bar{\theta}_j)\mu_j(\phi),
\end{equation*}
where $\mu_j$ and $m_j$ denote the Radon measure of the support of
$L_j$ and its multiplicity respectively.

Furthermore, the set $\{\bar{\theta}_1, \cdots, \bar{\theta}_N\}$ does
not depend on the sequence of rescalings chosen. }

\vspace{.1in}

\begin{corollary}\label{sing}
Let $(M,J,\bar{\omega},\bar{g},\Omega)$ be a compact almost
Calabi-Yau manifold of complex dimension $n$. Assume that the
initial surface $L_0$ is Lagrangian and zero-Maslov class. Then
there is no finite time Type-I singularity for generalized
Lagrangian mean curvature flow.
\end{corollary}

\vspace{.1in}

Note that, recently, Smoczyk-Wang (\cite{SW}) studied mean curvature
flow in even more general setting: generalized Lagrangian mean
curvature flow in symplectic manifold.

\vspace{.1in}

After we finished the first version of this paper, J. Li and the second author (\cite{LY}) considered the generalized symplectic mean curvature flow in
almost Einstein manifold. They proved that there is no Type I singularity for such a flow. They also showed global existence and convergence of the flow in graphic case.

\vspace{.1in}

\noindent \emph{Acknowledgement: The authors thank Professor Jiayu
Li for helpful discussions. The research was supported by NSFC No. 11071236. }

\vspace{.2in}

\section{Evolution Equations}

\vspace{.1in}

In this section, we will compute the evolution equations of the
induced metric and the second fundamental form of $L_t$ along the
generalized Lagrangian mean curvature flow (\ref{flow}).

\begin{lemma}\label{metric}
Along the generalized Lagrangian mean curvature flow (\ref{flow}),
the induced metric evolves by
\begin{equation*}
    \frac{\partial}{\partial t}g_{ij}=-2K^{\alpha}h^{\alpha}_{ij},
\end{equation*}
where $\textbf{K}=K^{\alpha}e_{\alpha}$.
\end{lemma}

{\it Proof.} We have
\begin{eqnarray*}
  \frac{\partial}{\partial
  t}g_{ij}&=&\langle\overline{\nabla}_{\textbf{K}}
  e_i,e_j\rangle+\langle e_i,\overline{\nabla}_{\textbf{K}} e_j\rangle
      =\langle\overline{\nabla}_{e_i}\textbf{K},e_j\rangle+\langle e_i,\overline{\nabla}_{e_j}\textbf{K}\rangle\\
  &=& K^\alpha\langle\overline{\nabla}_{e_i}e_\alpha,e_j\rangle+K^\alpha\langle
  e_i,\overline{\nabla}_{e_j}e_\alpha\rangle=-2K^\alpha h^\alpha_{ij}.
\end{eqnarray*}
This proves the lemma. \hfill Q.E.D.

\vspace{.1in}

As a corollary, we can immediately obtain:

\vspace{.1in}

\begin{corollary}
The area element of $L_t$ satisfies the following equation,
  \begin{equation}\label{e2.1}
    \frac{\partial}{\partial
    t}d\mu_t=-\langle\textbf{K},\textbf{H}\rangle d\mu_t,
  \end{equation}
  and consequently,
  \begin{equation}\label{e2.2}
    \frac{\partial}{\partial
    t}\int_{L_t}d\mu_t=-\int_{L_t}\langle\textbf{K},\textbf{H}\rangle
    d\mu_t.
  \end{equation}
\end{corollary}

Next, we compute the evolution equation of the second fundamental form. For the purpose of simplicity, we denote
\begin{equation}\label{e2.3}
    \textbf{K}=\textbf{H}-\textbf{V},
\end{equation}
where $\textbf{V}=n\pi_{\nu L}(\bar{\nabla}\psi)=V^{\alpha}e_{\alpha}$.

\begin{lemma}\label{lem2.3}
Along the generalized Lagrangian mean curvature flow (\ref{flow}), the second fundamental form
  $h^\alpha_{ij}$ satisfies
\begin{eqnarray}\label{h}
\frac{\partial}{\partial t}h^{\alpha}_{ij}&=&\Delta
h^\alpha_{ij}+(\overline{\nabla}_k \bar{R})_{\alpha
ijk}+(\overline{\nabla}_j \bar{R})_{\alpha
kik}-2\bar{R}_{lijk}h^\alpha_{lk}+2\bar{R}_{\alpha\beta
jk}h^\beta_{ik}+2\bar{R}_{\alpha\beta
ik}h^\beta_{jk}\nonumber\\
&&-\bar{R}_{lkik}h^\alpha_{lj}-\bar{R}_{lkjk}h^\alpha_{li}+\bar{R}_{\alpha
k\beta k}h^\beta_{ij}-h^\alpha_{im}(H^\gamma
h^\gamma_{mj}-h^\gamma_{mk}h^\gamma_{jk})\nonumber\\
&&-h^\alpha_{mk}(h^\gamma_{mj}h^\gamma_{ik}-h^\gamma_{mk}h^\gamma_{ij})
-h^\beta_{ik}(h^\beta_{lj}h^\alpha_{lk}-h^\beta_{lk}h^\alpha_{lj})-h^\alpha_{jk}h^\beta_{ik}H^\beta\nonumber\\&&-V^\alpha_{,ij}+V^\beta
h^\beta_{ik}h^\alpha_{jk}+V^\beta\bar{R}_{\beta
ji\alpha}+h^\beta_{ij}\langle e_\beta, \overline{\nabla}_{\textbf{K}}
e_\alpha\rangle,
\end{eqnarray}
where $\bar{R}_{ABCD}$ is the curvature tensor of $(M,\bar{g})$,
$\overline{\nabla}$ is the covariant derivative of $(M,\bar{g})$ and
$V^{\alpha}_{,ij}=\langle
\overline{\nabla}^N_{e_k}\overline{\nabla}^N_{e_i}\textbf{V},e_{\alpha}\rangle$.
Therefore, $|\textbf{A}|^2$ satisfies the following equation along the flow
(\ref{flow}),
\begin{eqnarray}\label{A}
\frac{\partial}{\partial t}|\textbf{A}|^{2}
    &= & \Delta |\textbf{A}|^{2} -2|\nabla \textbf{A}|^{2}
       +2[(\overline{\nabla}_k\bar{R})_{\alpha ijk}+(\overline{\nabla}_j\bar{R})_{\alpha kik}]h^{\alpha}_{ij} \nonumber\\
    & & -4\bar{R}_{lijk}h^{\alpha}_{lk}h^{\alpha}_{ij} +8\bar{R}_{\alpha\beta jk}h^{\beta}_{ik}h^{\alpha}_{ij}
       -4\bar{R}_{lkik}h^{\alpha}_{lj}h^{\alpha}_{ij}+ 2\bar{R}_{\alpha k\beta k}h^{\beta}_{ij}h^{\alpha}_{ij} \nonumber\\
    & & +2\sum_{\alpha,\beta,i,j}(\sum_{k}(h^{\alpha}_{ik}h^{\beta}_{jk}-h^{\alpha}_{jk}h^{\beta}_{ik}))^{2}
       +2\sum_{i,j,m,k}(\sum_{\alpha}h^{\alpha}_{ij}h^{\alpha}_{mk})^{2} \nonumber\\
    & & -2h^{\alpha}_{ij}V^{\alpha}_{,ij}+2V^{\beta}h^{\beta}_{ik}h^{\alpha}_{kj}h^{\alpha}_{ij}+2V^{\beta}h^{\alpha}_{ij}\bar{R}_{\beta ji\alpha}.
\end{eqnarray}
We also have,
\begin{eqnarray}\label{estimateA}
  \frac{\partial}{\partial t}|\textbf{A}|^2\leq\Delta|\textbf{A}|^2-|\nabla
  \textbf{A}|^2+C|\textbf{A}|^4+C|\textbf{A}|.
\end{eqnarray}
More generally, we have
\begin{eqnarray}\label{dA}
\frac{\partial}{\partial t}|\nabla^m \textbf{A}|^2\leq\Delta|\nabla^m
\textbf{A}|^2-|\nabla^{m+1}\textbf{A}|^2+C\sum_{i+j+k\leq m}|\nabla^i\textbf{A}||\nabla^j\textbf{A}||\nabla^k\textbf{A}||\nabla^m\textbf{A}|.
\end{eqnarray}
\end{lemma}
{\it Proof.} We will compute pointwise. So we choose normal
coordinate at a fixed point $p\in L$ such that
$\nabla_{e_j}e_i=(\overline{\nabla}_{e_j}e_i)^{T}=0$ at $p$. By
(7.4) in \cite{Wa}, the Laplacian of $h^\alpha_{ij}$ satisfies
\begin{eqnarray}\label{lh}
  \Delta h^\alpha_{ij}&=&H^\alpha_{,ij}-(\overline\nabla_k\bar{R})_{\alpha
  ijk}-(\overline{\nabla}_j\bar{R})_{\alpha ijk}+2\bar{R}_{lijk}h^\alpha_{lk}-2\bar{R}_{\alpha\beta
jk}h^\beta_{ik}-2\bar{R}_{\alpha\beta
ik}h^\beta_{jk}\nonumber\\
&&-\bar{R}_{\alpha
ij\beta}H^\beta+\bar{R}_{lkik}h^\alpha_{lj}+\bar{R}_{lkjk}h^\alpha_{li}-\bar{R}_{\alpha
k\beta k}h^\beta_{ij}+h^\alpha_{im}(H^\gamma
h^\gamma_{mj}-h^\gamma_{mk}h^\gamma_{jk})\nonumber\\
&&+h^\alpha_{mk}(h^\gamma_{mj}h^\gamma_{ik}-h^\gamma_{mk}h^\gamma_{ij})
+h^\beta_{ik}(h^\beta_{lj}h^\alpha_{lk}-h^\beta_{lk}h^\alpha_{lj}).
\end{eqnarray}
Now we compute $\frac{\partial}{\partial t}h^\alpha_{ij}.$ Since
$h^\alpha_{ij}=\langle\overline\nabla_{e_j}e_i,e_\alpha\rangle$, we have
\begin{eqnarray*}
\frac{\partial}{\partial
t}h^\alpha_{ij}&=&\langle\overline\nabla_{\textbf{K}}\overline\nabla_{e_j}e_i,e_\alpha\rangle
+\langle\overline{\nabla}_{e_j}e_i,\overline{\nabla}_{\textbf{K}}e_\alpha\rangle\\
&=&\langle\overline\nabla_{e_j}\overline{\nabla}_{\textbf{K}}e_i,e_\alpha\rangle-\langle\bar{R}(\textbf{K},e_j)e_i,e_\alpha\rangle
+\langle\overline{\nabla}_{e_j}e_i,\overline{\nabla}_{\textbf{K}}e_\alpha\rangle\\
&=&\langle\overline\nabla_{e_j}\overline{\nabla}_{e_i}\textbf{K},e_\alpha\rangle-\langle\bar{R}(\textbf{K},e_j)e_i,e_\alpha\rangle
+\langle\overline{\nabla}_{e_j}e_i,\overline{\nabla}_{\textbf{K}}e_\alpha\rangle.
\end{eqnarray*}
By breaking $\overline\nabla_{e_j}\overline\nabla_{e_i}\textbf{K}$
into normal and tangent parts, we get
\begin{eqnarray*}
  \langle\overline\nabla_{e_j}\overline\nabla_{e_i}\textbf{K},e_\alpha\rangle&=&
  \langle\overline\nabla_{e_j}[(\overline\nabla_{e_i}\textbf{K})^T+(\overline\nabla_{e_i}\textbf{K})^N],e_\alpha\rangle\\
  &=&\langle\overline\nabla_{e_j}^N\overline\nabla_{e_i}^N{\textbf{K}},e_\alpha\rangle-\langle(\overline\nabla_{e_i}\textbf{K})^T,\overline\nabla_{e_j}e_\alpha\rangle.
\end{eqnarray*}
Therefore,
\begin{eqnarray*}
  \frac{\partial}{\partial
  t}h^\alpha_{ij}=K^\alpha_{,ij}-K^\beta\bar{R}_{\beta
  ji\alpha}-\langle(\overline\nabla_{e_i}\textbf{K})^T,\overline\nabla_{e_j}e_\alpha\rangle
  +\langle\overline{\nabla}_{e_j}e_i,\overline{\nabla}_{\textbf{K}}e_\alpha\rangle,
\end{eqnarray*}
where $K^{\alpha}_{,ij}=\langle
\bar{\nabla}^N_{e_k}\bar{\nabla}^N_{e_i}\textbf{K},e_{\alpha}\rangle$. Note that
$\langle(\overline\nabla_{e_i}\textbf{K})^T,\overline\nabla_{e_j}e_\alpha\rangle=K^\beta
h^\beta_{ik}h^\alpha_{jk}$. By our choice of coordinate,
$\langle\overline{\nabla}_{e_j}e_i,\overline{\nabla}_{\textbf{K}}e_\alpha\rangle=h^\beta_{ij}\langle
e_\beta,\overline\nabla_{\textbf{K}}e_\alpha\rangle$. Therefore we have
\begin{eqnarray}\label{th}
  \frac{\partial}{\partial t}h^\alpha_{ij}=K^\alpha_{,ij}-K^\beta
  h^\beta_{ik}h^\alpha_{jk}-K^\beta\bar{R}_{\beta
  ji\alpha}+h^\beta_{ij}\langle
  e_\beta,\overline\nabla_{\textbf{K}}e_\alpha\rangle
\end{eqnarray}
Combine Equation (\ref{lh}) and (\ref{th}), we get the parabolic
equation (\ref{h}) for $h^\alpha_{ij}$.

\vspace{.1in}

Since $|\textbf{A}|^2=g^{ik}g^{jl}h^\alpha_{ij}h^\alpha_{kl}$, by (\ref{h})
we have,
\begin{eqnarray*}
  \frac{\partial}{\partial t}|\textbf{A}|^2&=&2(\frac{\partial}{\partial
  t}g^{ik})g^{jl}h^\alpha_{ij}h^\alpha_{kl}+2(\frac{\partial}{\partial
  t}h^\alpha_{ij})h^\alpha_{ij}\\
  &=&4(H^\beta-V^\beta)h^\beta_{ik}h^\alpha_{ij}h^\alpha_{kj}\\
  &&+2h^\alpha_{ij}[\Delta
h^\alpha_{ij}+(\overline{\nabla}_k \bar{R})_{\alpha
ijk}+(\overline{\nabla}_j \bar{R})_{\alpha
kik}-2\bar{R}_{lijk}h^\alpha_{lk}+2\bar{R}_{\alpha\beta
jk}h^\beta_{ik}+2\bar{R}_{\alpha\beta
ik}h^\beta_{jk}\nonumber\\
&&-\bar{R}_{lkik}h^\alpha_{lj}-\bar{R}_{lkjk}h^\alpha_{li}+\bar{R}_{\alpha
k\beta k}h^\beta_{ij}-h^\alpha_{im}(H^\gamma
h^\gamma_{mj}-h^\gamma_{mk}h^\gamma_{jk})\nonumber\\
&&-h^\alpha_{mk}(h^\gamma_{mj}h^\gamma_{ik}-h^\gamma_{mk}h^\gamma_{ij})
-h^\beta_{ik}(h^\beta_{lj}h^\alpha_{lk}-h^\beta_{lk}h^\alpha_{lj})-h^\alpha_{jk}h^\beta_{ik}H^\beta\nonumber\\&&-V^\alpha_{,ij}+V^\beta
h^\beta_{ik}h^\alpha_{jk}+V^\beta\bar{R}_{\beta
ji\alpha}+h^\beta_{ij}\langle e_\beta, \overline{\nabla}_{\textbf{K}}
e_\alpha\rangle].
\end{eqnarray*}
Using
\begin{eqnarray*}
  \Delta|\textbf{A}|^2=2|\nabla \textbf{A}|^2+2h^\alpha_{ij}\Delta h^\alpha_{ij}
\end{eqnarray*}
and the antisymmetric of $\langle e_\beta, \overline{\nabla}_{\textbf{K}}
e_\alpha\rangle$, and calculate similarly as in \cite{Wa}, we can
get the parabolic equation (\ref{A}) for $|\textbf{A}|^2$.

\vspace{.1in}

To prove (\ref{estimateA}), we only need to estimate the term
$h^\alpha_{ij}V^\alpha_{,ij}$. Note that $V^{\alpha}=n\langle\overline{\nabla}\psi, e_{\alpha}\rangle$. By definition, we know
\begin{eqnarray}\label{V1}
V^{\alpha}_{,i}
    & = & \langle \overline{\nabla}^N_{e_i}(V^{\gamma}e_{\gamma}),e_{\alpha}\rangle
          =e_i(V^{\gamma})\langle e_{\gamma},e_{\alpha}\rangle+V^{\gamma}\langle \overline{\nabla}^N_{e_i}e_{\gamma},e_{\alpha}\rangle\nonumber\\
    & = & e_i(V^{\alpha})+n\langle \overline{\nabla}\psi,e_{\gamma}\rangle\langle \overline{\nabla}^N_{e_i}e_{\gamma},e_{\alpha}\rangle
          \nonumber\\
    & = & n\overline{\nabla}_{e_i}\langle\overline{\nabla}\psi ,e_{\alpha}\rangle
             -n\langle \overline{\nabla}\psi,e_{\gamma}\rangle\langle e_{\gamma},\overline{\nabla}^N_{e_i}e_{\alpha}\rangle\nonumber\\
    & = & n\langle\overline{\nabla}_{e_i}\overline{\nabla}\psi ,e_{\alpha}\rangle+n\langle\overline{\nabla}\psi , \overline{\nabla}_{e_i}e_{\alpha}\rangle
             -n\langle \overline{\nabla}\psi,\overline{\nabla}^N_{e_i}e_{\alpha}\rangle\nonumber\\
    & = & n\langle\overline{\nabla}_{e_i}\overline{\nabla}\psi ,e_{\alpha}\rangle-n h^{\alpha}_{ik}\langle\overline{\nabla}\psi ,e_k\rangle.
\end{eqnarray}
Next, we compute the second covariant derivative at $p$. By
(\ref{V1}) and our choice of the frame:
\begin{eqnarray}\label{V2}
V^{\alpha}_{,ij}
    & = & \langle \overline{\nabla}^N_{e_j}\overline{\nabla}^N_{e_i}(V^{\gamma}e_{\gamma}),e_{\alpha}\rangle
                =\langle \overline{\nabla}^N_{e_j}(V^{\gamma}_{,i}e_{\gamma}),e_{\alpha}\rangle\nonumber\\
    & = & e_j(V^{\gamma}_{,i})\langle e_{\gamma},e_{\alpha}\rangle
            +V^{\gamma}_{,i}\langle \overline{\nabla}^N_{e_j}e_{\gamma},e_{\alpha}\rangle\nonumber\\
    & = & e_j(V^{\alpha}_{,i})
            +n\langle\overline{\nabla}_{e_i}\bar{\nabla}\psi ,e_{\gamma}\rangle\langle \overline{\nabla}^N_{e_j}e_{\gamma},e_{\alpha}\rangle
             -n h^{\gamma}_{ik}\langle\overline{\nabla}\psi ,e_k\rangle\langle \overline{\nabla}^N_{e_j}e_{\gamma},e_{\alpha}\rangle \nonumber\\
    & = & n\overline{\nabla}_{e_{j}}\langle\overline{\nabla}_{e_i}\overline{\nabla}\psi ,e_{\alpha}\rangle
          -ne_j(h^{\alpha}_{ik})\langle\overline{\nabla}\psi ,e_k\rangle
          -nh^{\alpha}_{ik}\overline{\nabla}_{e_{j}}\langle\overline{\nabla}\psi ,e_k\rangle\nonumber\\
    &   & -n\langle\overline{\nabla}_{e_i}\overline{\nabla}\psi ,e_{\gamma}\rangle\langle e_{\gamma},\overline{\nabla}^N_{e_j}e_{\alpha}\rangle
             -n h^{\gamma}_{ik}\langle\overline{\nabla}\psi ,e_k\rangle\langle \overline{\nabla}^N_{e_j}e_{\gamma},e_{\alpha}\rangle \nonumber\\
    & = & n\langle \overline{\nabla}_{e_{j}}\bar{\nabla}_{e_i}\overline{\nabla}\psi ,e_{\alpha}\rangle
          +n\langle\overline{\nabla}_{e_i}\overline{\nabla}\psi, \overline{\nabla}_{e_{j}}e_{\alpha}\rangle
          -ne_j(h^{\alpha}_{ik})\langle\overline{\nabla}\psi, e_k\rangle\nonumber\\
    &   & -nh^{\alpha}_{ik}\langle\overline{\nabla}_{e_{j}}\overline{\nabla}\psi ,e_k\rangle
          -nh^{\alpha}_{ik}\langle\overline{\nabla}\psi ,\overline{\nabla}_{e_{j}}e_k\rangle\nonumber\\
    &   & -n\langle\overline{\nabla}_{e_i}\overline{\nabla}\psi,\overline{\nabla}^N_{e_j}e_{\alpha}\rangle
             -n h^{\gamma}_{ik}\langle\overline{\nabla}\psi ,e_k\rangle\langle \overline{\nabla}^N_{e_j}e_{\gamma},e_{\alpha}\rangle\nonumber\\
    & = & -n(e_j(h^{\alpha}_{ik})+h^{\gamma}_{ik}\langle \overline{\nabla}^N_{e_j}e_{\gamma},e_{\alpha}\rangle)\langle\overline{\nabla}\psi ,e_k\rangle
          +n\langle\overline{\nabla}_{e_i}\overline{\nabla}\psi,\overline{\nabla}_{e_j}e_{\alpha}-\overline{\nabla}^N_{e_j}e_{\alpha}\rangle\nonumber\\
    &   & +n\langle\overline{\nabla}_{e_{j}}\overline{\nabla}_{e_i}\overline{\nabla}\psi ,e_{\alpha}\rangle
           -nh^{\alpha}_{ik}\langle\overline{\nabla}_{e_{j}}\overline{\nabla}\psi ,e_k\rangle
           -nh^{\alpha}_{ik}h^{\beta}_{kj}\langle\overline{\nabla}\psi,e_{\beta}\rangle\nonumber\\
    & = & -n(e_j(h^{\alpha}_{ik})+h^{\gamma}_{ik}\langle \overline{\nabla}^N_{e_j}e_{\gamma},e_{\alpha}\rangle)\langle\overline{\nabla}\psi ,e_k\rangle
          +n\langle\overline{\nabla}_{e_{j}}\overline{\nabla}_{e_i}\overline{\nabla}\psi ,e_{\alpha}\rangle\nonumber\\
    &   & -nh^{\alpha}_{jk}\langle\overline{\nabla}_{e_i}\overline{\nabla}\psi,e_k\rangle
           -nh^{\alpha}_{ik}\langle\overline{\nabla}_{e_{j}}\overline{\nabla}\psi ,e_k\rangle
           -nh^{\alpha}_{ik}h^{\beta}_{jk}\langle\overline{\nabla}\psi,e_{\beta}\rangle.
\end{eqnarray}
By the definition of covariant derivative of the second fundamental
form (Section 7 of \cite{Wa}) and the choice of frame, we have at
$p$
\begin{eqnarray}\label{V3}
  h^{\alpha}_{ik,j} &=& \langle (\overline{\nabla}_{e_j}\textbf{A})(e_i,e_k),e_{\alpha}\rangle \nonumber\\
   &=&  \langle \overline{\nabla}_{e_j}^{N}(\textbf{A}(e_i,e_k))-\textbf{A}((\overline{\nabla}_{e_j}e_i)^T,e_k)
               -\textbf{A}(e_i,(\overline{\nabla}_{e_j}e_k)^T),e_{\alpha}\rangle\nonumber\\
   &=& \langle \overline{\nabla}_{e_j}^{N}(h^{\gamma}_{ik}e_{\gamma}),e_{\alpha}\rangle
        =e_{j}(h^{\alpha}_{ik})+h^{\gamma}_{ik}\langle \overline{\nabla}^N_{e_j}e_{\gamma},e_{\alpha}\rangle.
\end{eqnarray}
Therefore, we have
\begin{eqnarray*}
V^{\alpha}_{,ij}
    & = & -nh^{\alpha}_{ik,j}\langle\overline{\nabla}\psi ,e_k\rangle
          +n\langle\overline{\nabla}_{e_{j}}\overline{\nabla}_{e_i}\overline{\nabla}\psi ,e_{\alpha}\rangle
          -nh^{\alpha}_{jk}\langle\overline{\nabla}_{e_i}\overline{\nabla}\psi,e_k\rangle\nonumber\\
    &   & -nh^{\alpha}_{ik}\langle\overline{\nabla}_{e_{j}}\overline{\nabla}\psi ,e_k\rangle
           -nh^{\alpha}_{ik}h^{\beta}_{kj}\langle\overline{\nabla}\psi,e_{\beta}\rangle,
\end{eqnarray*}
and
\begin{eqnarray*}
h^{\alpha}_{ij}V^{\alpha}_{,ij}
    & = & -nh^{\alpha}_{ij}h^{\alpha}_{ik,j}\langle\overline{\nabla}\psi ,e_k\rangle
          +nh^{\alpha}_{ij}\langle\overline{\nabla}_{e_{j}}\overline{\nabla}_{e_i}\overline{\nabla}\psi ,e_{\alpha}\rangle
          -2nh^{\alpha}_{ij}h^{\alpha}_{jk}\langle\overline{\nabla}_{e_i}\overline{\nabla}\psi,e_k\rangle\nonumber\\
    &   & -nh^{\alpha}_{ij}h^{\alpha}_{ik}h^{\beta}_{kj}\langle\overline{\nabla}\psi,e_{\beta}\rangle.
\end{eqnarray*}
Note that
\begin{eqnarray*}
   |\langle\overline{\nabla}\psi ,e_j\rangle|\leq ||\psi||_{C^1(M)},
    & & |\langle\overline{\nabla}_{e_{k}}\overline{\nabla}_{e_i}\overline{\nabla}\psi ,e_{\alpha}\rangle|\leq ||\psi||_{C^3(M)}, \\
   |\langle\overline{\nabla}_{e_i}\overline{\nabla}\psi,e_j\rangle| \leq ||\psi||_{C^2(M)},
   & &  |\langle\overline{\nabla}\psi,e_{\beta}\rangle|\leq ||\psi||_{C^1(M)}.
\end{eqnarray*}
Thus, by Young's inequality, we have
\begin{equation}\label{e2.13}
    |h^{\alpha}_{ij}V^{\alpha}_{,ij}|\leq \varepsilon|\nabla \textbf{A}|^2+C_1(\varepsilon)|\textbf{A}|^2
     +C_2|\textbf{A}|+C_3|\textbf{A}|^3.
\end{equation}
As $\textbf{V}$ is bounded, combining (\ref{A}) with (\ref{e2.13}) yields (\ref{estimateA}).

\vspace{.1in}

For the higher derivative estimate, we only need to notice that by induction, we can show that for each $m$
\begin{eqnarray*}
 \frac{\partial}{\partial t}\nabla^{m}\textbf{A}
    & = & \Delta\nabla^{m}\textbf{A}+\sum_{i+j+k=m}\nabla^{i}\textbf{A}*\nabla^{j}\textbf{A}*\nabla^{k}\textbf{A}
           +\sum_{i+j+k=m}\nabla^{i}\textbf{A}*\nabla^{j}\textbf{A}*\nabla^{k}\textbf{V}\nonumber\\
    &   & +\sum^{m}_{j=0}\sum_{\mbox{\tiny %
      $ \begin{array}{c}
          i_{1}+\cdots+i_{r}+k=m+1-j\\
          i_{1},\cdots,i_{r},k\geq 1
        \end{array}
      $ } }\overline{\nabla}^{j}\bar{R}*\nabla^{i_{1}-1}\textbf{A}*\cdots *\nabla^{i_{r}-1}\textbf{A}*\nabla^{k-1}\textbf{V}\nonumber\\
    &   & + \nabla ^{m+2}\textbf{V}+\overline{\nabla}^{m+1}\bar{R}.
\end{eqnarray*}
and
\begin{eqnarray}\label{V}
 |\nabla^{m}\textbf{V}|
    & \leq & C\sum_{\mbox{\tiny %
      $ \begin{array}{c}
          i_{1}+\cdots+i_{r}+k=m+1\\
          i_{1},\cdots,i_{r},k\geq 1
        \end{array}
      $ } }|\nabla^{i_{1}-1}\textbf{A}| \cdots |\nabla^{i_{r}-1}\textbf{A}| |\overline{\nabla}^{k}\psi|.
\end{eqnarray}
Then arguing in the same way as that of the mean curvature flow (see, for example, Section 3 of \cite{HS}), we can obtain (\ref{dA}).
\hfill Q.E.D

\vspace{.1in}

Once we have Lemma \ref{lem2.3}, we can get the following longtime
existence theorem. This is essentially the same as in mean curvature
flow case (see for example, Lemma 7.2 of \cite{Hu3}).

\begin{theorem}
If the second fundamental form of $L_t$ is uniformly bounded under
  the generalized Lagrangian mean curvature flow (\ref{flow}) for
  all time $t\in[0,T)$, then the solution can be extended beyond
  $T$.
\end{theorem}

{\it Proof.} If $|\textbf{A}|(t)\leq C$ for $t\in [0,T)$, then by
(\ref{dA}) and the standard application of parabolic maximum
principle, we know that
\begin{equation}\label{e2.15}
    |\nabla^m\textbf{A}|\leq C(m), \ \ for \ t\in [0,T)
\end{equation}
for some constant $C(m)$. Then (\ref{V}) implies that
\begin{equation}\label{e2.16}
    |\nabla^m\textbf{V}|\leq C(m), \ \ for \ t\in [0,T).
\end{equation}
Therefore, by definition,
\begin{equation}\label{e2.17}
    |\nabla^m\textbf{K}|\leq C(m), \ \ for \ t\in [0,T).
\end{equation}
By the equation (\ref{flow}), the remaining part of the proof is standard and we omit the details here.
\hfill Q.E.D.

\vspace{.1in}

Using (\ref{estimateA}), we can argue in the same way as in the mean curvature flow (for example, Lemma 4.6 of \cite{CL1}) to obtain the lower bound of the blow up rate of the maximal norm of the second fundamental form at finite singular time $T$:

\begin{lemma}
Let $U_{t}=\max_{M_{t}}|A|^{2}$. If the generalized Lagrangian mean
curvature flow (\ref{flow}) blows up at finite time $T>0$, there is
a positive $c$ depending only on $M$, such that if $0<T-t<
\frac{\pi}{32\sqrt{c}}$, then the function $U_{t}$ satisfies
\begin{equation*}
U_{t}\geq \frac{1}{8\sqrt{2}(T-t)}.
\end{equation*}
\end{lemma}

According to the upper bound of the blow up rate, we can classify
the singularities of the generalized Lagrangian mean curvature flow
(\ref{flow}) into two types, which is similar to that of mean
curvature flow defined by Huisken (\cite{Hu2}).

\begin{definition}
We say that the generalized Lagrangian mean curvature flow
(\ref{flow}) develops Type-I singularity at $T>0$, if
\begin{equation*}
\limsup_{t\to T}(T-t)\max_{M_{t}}|A|^{2}\leq C,
\end{equation*}
for some positive constant C. Otherwise, we say the generalized
Lagrangian mean curvature flow (\ref{flow}) develops Type-II
singularity.
\end{definition}

\vspace{.2in}

\section{A Monotonicity Formula}

\vspace{.1in}

Let $H(\textbf{X},\textbf{X}_0,t_0,t)$ be the backward heat kernel
on ${\textbf{R}}^{k}$. Let $L_t$ be a smooth family of submanifolds
of dimension $n$ in ${\textbf{R}}^k$ defined by $F_t: L\rightarrow
{\textbf{R}}^k$. Define
\begin{equation*}
    \rho(\textbf{X},t)=(4\pi(t_0-t))^{\frac{k-n}{2}}H(\textbf{X},\textbf{X}_0,t_0,t)
  =\frac{1}{(4\pi(t_0-t))^{\frac{n}{2}}}\exp\left(-\frac{|\textbf{X}-\textbf{X}_0|^2}{4(t_0-t)}\right)
\end{equation*}
for $t<t_0$. We have along the generalized mean curvature flow
(\ref{flow})
\begin{eqnarray}\label{e3.1}
      \frac{\partial\rho}{\partial t}
  &=&-\frac{n}{2}\frac{-4\pi}{(4\pi(t_0-t))^{\frac{n}{2}+1}}\exp\left(-\frac{|\textbf{X}-\textbf{X}_0|^2}{4(t_0-t)}\right)
        +\rho(x,t)\frac{\langle \textbf{K},\textbf{X}-\textbf{X}_0\rangle}{-2(t_0-t)}
           +\rho(x,t)\frac{-|\textbf{X}-\textbf{X}_0|^2}{4(t_0-t)^2}\nonumber\\
  &=&\left(\frac{n}{2(t_0-t)}-\frac{\langle \textbf{K}, \textbf{X}-\textbf{X}_0\rangle}{2(t_0-t)}-\frac{|\textbf{X}-\textbf{X}_0|^2}{4(t_0-t)^2}\right)\rho.
\end{eqnarray}
As
\begin{equation*}
  \nabla\exp\left(-\frac{|\textbf{X}-\textbf{X}_0|^2}{4(t_0-t)}\right)
  =-\exp\left(-\frac{|\textbf{X}-\textbf{X}_0|^2}{4(t_0-t)}\right)\frac{\langle\textbf{X}-\textbf{X}_0,\nabla\textbf{X}\rangle}{2(t_0-t)},
\end{equation*}
we have
\begin{equation*}
  \Delta\rho=\left(\frac{\langle\textbf{X}-\textbf{X}_0,\nabla\textbf{X}\rangle^2}{4(t_0-t)^2}
         -\frac{\langle\textbf{X}-\textbf{X}_0,\Delta\textbf{X}\rangle}{2(t_0-t)}
      -\frac{|\nabla\textbf{X}|^2}{2(t_0-t)}\right)\rho.
\end{equation*}
Note that
\begin{equation*}
    |\nabla F|^2=n, \ \ \Delta F=\textbf{H}+g^{ij}\overline\Gamma^\alpha_{\rho\sigma}\frac{\partial
F^\rho}{\partial x^i}\frac{\partial F^\sigma}{\partial x^j}e_\alpha,
\end{equation*}
where $e_\alpha$, $\alpha=1,\cdots,n$ is a basis of
$T^\perp L_t$, $g^{ij}$ is the inverse of the induced metric on $L_t$ and
$\overline{\Gamma}^\alpha_{\rho\sigma}$ is the Christoffel symbol on $M$. Therefore, we
have
\begin{equation}\label{e3.2}
  \Delta\rho=\left(\frac{\langle F-\textbf{X}_0,\nabla F\rangle^2}{4(t_0-t)^2}-\frac{\langle F-\textbf{X}_0,\textbf{H}+g^{ij}\overline\Gamma^\alpha_{\rho\sigma}\frac{\partial
F^\rho}{\partial x^i}\frac{\partial F^\sigma}{\partial x^j}e_\alpha\rangle}{2(t_0-t)}
      -\frac{n}{2(t_0-t)}\right)\rho.
\end{equation}
Combining (\ref{e3.1}) with (\ref{e3.2}) gives us
\begin{eqnarray}\label{e3.3}
\left(\frac{\partial}{\partial t}+\Delta\right)\rho
  & = & \left(-\left|\textbf{K}+\frac{(F-\textbf{X}_0)^\perp}{2(t_0-t)}\right|^2+|\textbf{K}|^2-\frac{\langle
  g^{ij}\overline\Gamma^\alpha_{\rho\sigma}\frac{\partial
F^\rho}{\partial x^i}\frac{\partial F^\sigma}{\partial x^j}e_\alpha,F-\textbf{X}_0\rangle}{t_0-t} \right.\nonumber \\
  &  & \left. -\frac{n\langle\pi_{\nu L}(\overline{\nabla}\psi),F-\textbf{X}_0\rangle}{2(t_0-t)}\right)\rho.
\end{eqnarray}

\vspace{.1in}

Denote the injectivity radius of
$(M,\bar{g})$ by $i_M$. For $\textbf{X}_0\in M$, take a normal coordinate
neighborhood $U$ and let $\phi\in C_0^\infty(B_{2r}(\textbf{X}_0))$
be a cut-off function with $\phi\equiv1$ in $B_r(\textbf{X}_0)$,
$0<2r<i_M$. Using the local coordinates in $U$ we may regard
$F(x,t)$ as a point in ${\textbf{R}}^{k}$ whenever $F(x,t)$ lies in
$U$.

\vspace{.1in}

The following monotonicity formula generalizes Proposition 2.1 of
\cite{CL2} to the almost Calabi-Yau case.

\begin{proposition}\label{prop3.1}
  Let $F_t: L\rightarrow M$ be a smooth mean curvature flow of
  a compact Lagrangian submanifold $L_0$ in a compact almost
  Calabi-Yau manifold $M$ of complex dimension $n$. Let $f_t:L_t\to M$ be a family of smooth function for $t\in [0,T)$ which is uniformly bounded, then for any $\varepsilon>0$, there are
  positive constants $c_1$, $c_2$, $c_3$ and $c_4$ depending only on $M$,$F_0$, $t_0$, $\sup_{0\leq t<T}||f_t||_{C^{0}(L_t)}$, $\varepsilon$ and $r$ which is the constant in the definition of $\Psi$, such that
  \begin{eqnarray}\label{mono}
     \frac{\partial}{\partial t}\left(e^{c_1\sqrt{t_0-t}}\int_{L_t}f_t\phi\rho d\mu_t\right)
    &\leq&e^{{c_1}\sqrt{t_0-t}}\int_{L_t}(\frac{\partial}{\partial t}-\Delta)f_t\phi\rho d\mu_t\nonumber\\
    &&-e^{c_1\sqrt{t_0-t}}\int_{L_t}f_t\phi\rho
    \left|\textbf{K}+\frac{(F-\textbf{X}_0)^\perp}{2(t_0-t)}\right|^2d\mu_t\nonumber\\
    &&+\varepsilon
    e^{c_1\sqrt{t_0-t}}\int_{L_t}f_t\phi\rho|\textbf{K}|^2d\mu_t\nonumber\\
    &&+\frac{c_2 e^{c_1\sqrt{t_0-t}}}{(t_0-t)^{\frac{3}{4}}}+c_3 e^{c_1\sqrt{t_0-t}}.
  \end{eqnarray}
\end{proposition}

\vspace{.1in}

{\it Proof.} We define
\begin{eqnarray*}
  \Phi(\textbf{X}_0,t_0,t)=\int_{L_t}f_t\phi(F)\rho(F,\textbf{X}_0
  ,t,t_0)d\mu_t
\end{eqnarray*}
where $\rho$ is defined as above by taking $k=2n$. Note that
\begin{eqnarray*}
  \frac{\partial\phi(F)}{\partial
  t}=\overline\nabla\phi\cdot\textbf{K}.
\end{eqnarray*}
Using (\ref{e2.1}), (\ref{e2.3}) and (\ref{e3.3}) , we have
\begin{eqnarray}\label{e3.5}
  &&\frac{d}{dt}\int_{L_t}f_t\phi\rho d\mu_t\nonumber\\
  &=&\int_{L_t}\frac{\partial}{\partial
  t}f_t\phi\rho+\int_{L_t}f_t\langle\overline{\nabla}\phi,\textbf{K}\rangle\rho+\int_{L_t}f_t\phi\frac{\partial}{\partial
  t}\rho+\int_{L_t}f_t\phi\rho\frac{\partial}{\partial
  t}d\mu_t\nonumber\\
  &=&\int_{L_t}(\frac{\partial}{\partial
  t}-\Delta)f_t\phi\rho+\int_{L_t}\Delta f_t\phi\rho+\int_{L_t}f_t\langle\overline{\nabla}\phi,\textbf{K}\rangle\rho+\int_{L_t}f_t\phi(\frac{\partial}{\partial
  t}+\Delta)\rho-\int_{L_t}f_t\phi\Delta\rho\nonumber\\
  &&-\int_{L_t}f_t\phi\rho\langle\textbf{H},\textbf{K}\rangle\nonumber\\
  &\leq&\int_{L_t}(\frac{\partial}{\partial
  t}-\Delta)f_t\phi\rho+\int_{L_t}(\phi\rho\Delta
  f_t-f_t\phi\Delta\rho)+\int_{L_t}f_t\rho(\frac{\varepsilon}{2}\phi|\textbf{K}|^2+\frac{1}{2\varepsilon}\frac{|\overline\nabla\phi|^2}{\phi})\nonumber\\
  &&+\int_{L_t}\left(-\left|\textbf{K}+\frac{(F-\textbf{X}_0)^\perp}{2(t_0-t)}\right|^2+|\textbf{K}|^2
  -\frac{\langle g^{ij}\overline\Gamma^\alpha_{\rho\sigma}\frac{\partial
F^\rho}{\partial x^i}\frac{\partial F^\sigma}{\partial x^j}e_\alpha,
  F-\textbf{X}_0\rangle}{t_0-t}-\frac{\langle\textbf{V},F-\textbf{X}_0\rangle}{2(t_0-t)}\right)f_t\phi\rho\nonumber\\
  &&              -\int_{L_t}f_t\phi\rho\langle\textbf{K},\textbf{K}+\textbf{V}\rangle\nonumber\\
  &\leq&\int_{L_t}(\frac{\partial}{\partial
  t}-\Delta)f_t\phi\rho-\int_{L_t}f_t\phi\rho\left|\textbf{K}+\frac{(F-\textbf{X}_0)^\perp}{2(t_0-t)}\right|^2+\int_{L_t}(\phi\rho\Delta
  f_t-f_t\phi\Delta\rho)\nonumber\\
  &&-\int_{L_t}f_t\phi\rho\frac{\langle
  g^{ij}\overline\Gamma^\alpha_{\rho\sigma}\frac{\partial
F^\rho}{\partial x^i}\frac{\partial F^\sigma}{\partial
x^j}e_\alpha,F-\textbf{X}_0\rangle}{t_0-t}+\frac{\varepsilon}{2}\int_{L_t}f_t\phi\rho|\textbf{K}|^2
+\frac{1}{2\varepsilon}\int_{L_t}f_t\rho\frac{|\overline\nabla\phi|^2}{\phi}\nonumber\\&&
  -\int_{L_t}f_t\phi\rho\frac{\langle\textbf{V},F-\textbf{X}_0\rangle}{2(t_0-t)}
  -\int_{L_t}f_t\phi\rho\langle\textbf{K},\textbf{V}\rangle.
\end{eqnarray}
Again, by (\ref{e2.1}) and (\ref{e2.3}), we have
\begin{eqnarray*}
  \frac{\partial}{\partial
  t}d\mu_t&=&-\langle\textbf{K},\textbf{H}\rangle
  d\mu_t=-\langle\textbf{H},\textbf{H}-\textbf{V}\rangle d\mu_t=(-|\textbf{H}|^2+\langle\textbf{H},\textbf{V}\rangle)d\mu_t
  \leq\frac{1}{4}|\textbf{V}|^2d\mu_t\leq
  Cd\mu_t,
\end{eqnarray*}
which implies that
\begin{eqnarray*}
  \frac{\partial}{\partial t}Area(L_t)\leq C Area(L_t).
\end{eqnarray*}
Therefore, we have
\begin{eqnarray}\label{area}
  Area(L_t)\leq e^{Ct_0}Area(L_0)\leq C.
\end{eqnarray}
By Stokes' theorem,
\begin{equation*}
    \int_{L_t}\left(\phi\rho\Delta f_t-f_t\phi\Delta\rho\right)
    =\int_{L_t}f_t\Delta\phi\rho(F,t)+2\int_{L_t}f_t\langle\nabla\phi,\nabla\rho(F,t)\rangle.
\end{equation*}
Note that $\Delta\phi=0$, $\nabla\phi=0$ in $B_r(\textbf{X}_0)$, we
can see that $\Delta\phi\rho(F,t)\leq C$ and
$|\langle\nabla\phi,\nabla\rho(F,t)\rangle|\leq C$. Hence
\begin{eqnarray*}
  \int_{L_t}f_t\Delta\phi\rho(F,t)d\mu_t\leq
  C\int_{L_t}d\mu_t\leq C,
\end{eqnarray*}
\begin{eqnarray*}
  \int_{L_t}f_t\langle\nabla\phi,\nabla\rho(F,t)\rangle d\mu_t\leq
  C\int_{L_t}d\mu_t\leq C.
\end{eqnarray*}
Therefore, we have
\begin{eqnarray}\label{e3.7}
\int_{L_t}\left(\phi\rho\Delta f_t-f_t\phi\Delta\rho\right)\leq C.
\end{eqnarray}
As $\phi\in C_0^\infty(B_{2r}(\textbf{X}_0), {\textbf{R}}^+)$, we have (Lemma 6.6 of \cite{Il})
\begin{equation}\label{e3.8}
    \frac{|\overline\nabla\phi|^2}{\phi}\leq 2\max_{\phi>0}|\overline{\nabla}^2 \phi|.
\end{equation}
By Young's inequality,
\begin{eqnarray}\label{e3.9}
  -\int_{L_t}f_t\phi\rho\langle\textbf{K},\textbf{V}\rangle\leq\frac{\varepsilon}{2}\int_{L_t}f_t\phi\rho|\textbf{K}|^2
  +C(\varepsilon)\int_{L_t}f_t\phi\rho|\textbf{V}|^2\leq\frac{\varepsilon}{2}\int_{L_t}f_t\phi\rho|\textbf{K}|^2+C(\varepsilon).
\end{eqnarray}
Since we choose a normal coordinates in
$B_{2\tilde{r}}(\textbf{X}_0)$ in $(M,\bar{g}(t))$, we have
$\overline\Gamma^\alpha_{\rho\sigma}(\textbf{X}_0,t)=0$, and
$|g^{ij}\overline\Gamma^\alpha_{\rho\sigma}\frac{\partial
F^\rho}{\partial x^i}\frac{\partial F^\sigma}{\partial x^j}|\leq
C|F-\textbf{X}_0|$, thus
\begin{eqnarray*}
\frac{\langle
F-\textbf{X}_0,g^{ij}\overline\Gamma^\alpha_{\rho\sigma}\frac{\partial
F^\rho}{\partial x^i}\frac{\partial F^\sigma}{\partial
x^j}e_\alpha\rangle}{2(t_0-t)}\leq
C\frac{|F-\textbf{X}_0|^2}{2(t_0-t)}.
\end{eqnarray*}
Similar to the proof of (13) in \cite{CL2}, we have
\begin{eqnarray}\label{e3.11}
\frac{\langle
F-\textbf{X}_0,g^{ij}\overline\Gamma^\alpha_{\rho\sigma}\frac{\partial
F^\rho}{\partial x^i}\frac{\partial F^\sigma}{\partial
x^j}e_\alpha\rangle}{2(t_0-t)}\rho(F,t)\leq
C_1\frac{\rho(F,t)}{\sqrt{t_0-t}}+C.
\end{eqnarray}
Finally, we need to estimate the term
$-\int_{L_t}f_t\phi\rho\frac{\langle\textbf{V},F-\textbf{X}_0\rangle}{2(t_0-t)}$.
We claim
\begin{eqnarray}\label{e3.12}
  \frac{|F-\textbf{X}_0|^2}{(t_0-t)^{\alpha}}\rho(F,t)\leq
  C_1\frac{\rho(F,t)}{(t_0-t)^{\beta}}+C, \ \ 0<\alpha-1 <\beta <1.
\end{eqnarray}
In fact it suffices to show for any $x$ and $s>0$
\begin{eqnarray*}
  \frac{x^2}{s^\alpha}\frac{e^{-x^2/s}}{s^{n/2}}\leq
  C\left(1+\frac{1}{s^\beta}\frac{e^{-x^2/s}}{s^{n/2}}\right).
\end{eqnarray*}
Let $y=x^2/s$ and then it suffices to show
\begin{eqnarray*}
  \frac{y}{s^{\alpha-1}}\frac{e^{-y}}{s^{n/2}}\leq
  C\left(1+\frac{1}{s^\beta}\frac{e^{-y}}{s^{n/2}}\right),
\end{eqnarray*}
which is equivalent to
\begin{eqnarray*}
  y\leq C\left(s^{n/2+\alpha-1}e^y+\frac{1}{s^{\beta+1-\alpha}}\right).
\end{eqnarray*}
If $y\leq\frac{1}{s^{\beta+1-\alpha}}$, then it hold trivially. If $y>\frac{1}{s^{\beta+1-\alpha}}$, then from $y^\gamma\leq
C(\gamma)e^y$ ($\gamma>1$ is to be determined), we see that
$y\leq\frac{C}{y^{\gamma-1}}e^y\leq Cs^{(\beta+1-\alpha)(\gamma-1)}e^y$. We only need to choose $\gamma$ such that
$(\beta+1-\alpha)(\gamma-1)=\frac{n}{2}+\alpha-1$, i.e.,
$\gamma=\frac{\frac{n}{2}+\beta}{\beta+1-\alpha}>\frac{\frac{n}{2}}{1+1-1}=\frac{n}{2}\geq1$.

As $\textbf{V}$ is bounded, we have
\begin{eqnarray*}
  -\int_{L_t}f_t\phi\rho\frac{\langle\textbf{V},F-\textbf{X}_0\rangle}{2(t_0-t)}
  & \leq &
  C\int_{L_t}f_t\phi\frac{|F-\textbf{X}_0|}{t_0-t}\rho\nonumber \\
  &\leq& C\left(\int_{L_t}\frac{|F-\textbf{X}_0|^2}{(t_0-t)^\alpha}f_t\phi\rho\right)^
  {\frac{1}{2}}\left(\int_{L_t}\frac{f_t\phi\rho}{(t_0-t)^{2-\alpha}}\right)^{\frac{1}{2}}\nonumber\\
  &\leq&C\int_{L_t}\frac{|F-\textbf{X}_0|^2}{(t_0-t)^\alpha}f_t\phi\rho+\frac{C}{(t_0-t)^{2-\alpha}}\nonumber\\
  &\leq&C\int_{L_t}f_t\phi\frac{\rho}{(t_0-t)^\beta}+\frac{C}{(t_0-t)^{2-\alpha}}+C\nonumber\\
  &\leq&\frac{C}{(t_0-t)^\beta}\Phi+\frac{C}{(t_0-t)^{2-\alpha}}+C.
\end{eqnarray*}
Especially, if we choose $\alpha=\frac{5}{4}$, $\beta=\frac{1}{2}$, then we have
\begin{eqnarray}\label{e3.13}
  -\int_{L_t}f_t\phi\rho\frac{\langle\textbf{V},F-\textbf{X}_0\rangle}{2(t_0-t)}
  \leq\frac{C}{\sqrt{t_0-t}}\Phi+\frac{C}{(t_0-t)^{\frac{3}{4}}}.
\end{eqnarray}
Putting (\ref{e3.7}), (\ref{e3.8}), (\ref{e3.9}), (\ref{e3.11}) and
(\ref{e3.13}) into (\ref{e3.5}), we obtain
\begin{eqnarray}\label{e3.14}
  \frac{\partial}{\partial t}\Phi
   & \leq & \int_{L_t}(\frac{\partial}{\partial t}-\Delta)f_t\phi\rho-\int_{L_t}f_t\phi\rho\left|\textbf{K}
   +\frac{(F-\textbf{X}_0)^\perp}{2(t_0-t)}\right|^2+\varepsilon\int_{L_t}f_t\phi\rho|\textbf{K}|^2\nonumber\\
   &      &  +\frac{c_1}{\sqrt{t_0-t}}\Phi+\frac{c_2}{(t_0-t)^{\frac{3}{4}}}+c_3.
\end{eqnarray}
Rearranging (\ref{e3.14}) yields the desired inequality.
\hfill Q.E.D.

\vspace{.2in}

\section{Proof of the Main Theorem}

\vspace{.1in}

Let $(L_t)_{0\leq t<T}$ be a smooth solution of the generalized Lagrangian mean curvature in an almost Calabi-Yau manifold with zero-Maslov class.  Recall that the rescaled flow is defined by
\begin{equation}\label{e3.15}
    F_i(x,t)\equiv F_{\lambda_i}(x,t)=\lambda_i(F(x,T+\lambda_i^{-2}s)-X_0), \ \ \ for \ -\lambda^2 T<s<0.
\end{equation}
Denote by $L^k_s$ the scaled surface $F_k(\cdot,s)$, then the induced metric satisfies
\begin{equation*}
g^i_{kl}=\lambda_i^2g_{kl},\ \ (g^i)^{kl}=\lambda_i^{-2}g^{kl}.
\end{equation*}
Moreover, it is easy to show that the scaled
surface also evolves by a generalized mean curvature flow
\begin{equation*}
\frac{\partial F_i}{\partial s}=\textbf{K}_i,
\end{equation*}
where
\begin{equation}\label{e3.16}
    \textbf{K}_i=\textbf{H}_i-\lambda_i^{-1}n\pi_{\nu L^i}(\overline{\nabla}\psi).
\end{equation}
Note that the Lagrangian angle $\theta_i$ on the rescaled surface $L^i_s$ satisfies
\begin{equation*}
    \theta_i(F_i(x,s))=\theta(F(x, T+\lambda_i^{-2}s)).
\end{equation*}

Proposition 5.1 in \cite{N} can be easily generalized to our case
that the ambient space is an almost Calabi-Yau manifold. The proof
is the same as in \cite{N}, so we just state the result here without proof.
\begin{proposition}\label{Li}
  Let $(L^i)$ be a sequence of smooth zero-Maslov class Lagrangians
  in $M$ such that, for some fixed $R>0$, the following properties
  hold:

  \noindent(a) There exists a constant $D_0$ for which
  $$\mu(L^i\cap B_{2R}(0))\leq D_0R^n  \ \ and \ \ \sup_{L^i\cap B_{2R}(0)}|\theta_i|\leq D_0$$
  for all $i\in{\textbf{N}}$.

  \noindent(b)
  \[\lim_{i\rightarrow\infty}\mu(\partial L^i\cap B_{2R}(0))=0\]
  and
  \[\lim_{i\rightarrow\infty}\int_{L^i\cap B_{2R}(0)}|\textbf{H}|^2d\mu=0.\]
Then there exist a finite set $\bar{\theta}_1, \cdots,
\bar{\theta}_N$ and integral Special Lagrangians
$$L_1, ..., L_N$$
such that, after passing to a subsequence, we have for every smooth
function $\phi$ compactly supported, every $f$ in $C^2(\textbf{R})$,
and every $s<0$
\begin{equation*}
    \lim_{i\to\infty}\int_{L^{i}_{s}}f(\theta_{i,s})\phi d\mu^{i}_{s}=\sum_{j=1}^{N}m_jf(\bar{\theta}_j)\mu_j(\phi),
\end{equation*}
where $\mu_j$ and $m_j$ denote the Radon measure of the support of
$L_j$ and its multiplicity respectively.
\end{proposition}

To prove the Main Theorem, we need the following lemma which generalizes Lemma 5.4 of \cite{N}:

\begin{lemma}\label{lem3.3}
  For any $s_1<s_2<0$ and for any $R>0$, we have
  \begin{eqnarray}\label{limit}
    \lim_{i\rightarrow\infty}\int_{s_1}^{s_2}\int_{L_s^i\cap
    B_R(0)}(|F_i^{\perp}|^2+|\textbf{H}_i|^2+|\textbf{K}_i|^2)d\mu_s^ids=0.
  \end{eqnarray}
\end{lemma}

{\it Proof.} As both $L$ and $M$ are compact, $\theta_0$ is bounded. Applying parabolic maximum principle to (\ref{e1.3}) yields that $\theta_t$ is uniformly bounded at each time. In fact, it is can be bounded in term of the bound of $\theta_0$.
By (\ref{e1.3}), we have
\begin{equation*}
    \frac{\partial}{\partial t}\theta^2 =\Delta(\theta^2)-2|\nabla\theta|^2+2n\theta d\psi(\nabla\theta)
      =\Delta(\theta^2)-2|\nabla\theta|^2+2n\theta \langle\nabla\psi,\nabla\theta\rangle.
\end{equation*}
Take $f_t=\theta_t^2$ in (\ref{mono}), we have
  \begin{eqnarray}\label{monotheta}
     \frac{\partial}{\partial t}\left(e^{c_1\sqrt{T-t}}\int_{L_t}\theta^2\phi\rho d\mu_t\right)
    &\leq&e^{{c_1}\sqrt{T-t}}\int_{L_t}(-2|\nabla\theta|^2\phi\rho+2n\theta\langle\nabla\psi,\nabla\theta\rangle\phi\rho) d\mu_t\nonumber\\
    &&-e^{c_1\sqrt{T-t}}\int_{L_t}\theta^2\phi\rho
    \left|\textbf{K}+\frac{(F-\textbf{X}_0)^\perp}{2(T-t)}\right|^2d\mu_t\nonumber\\
    &&+\varepsilon
    e^{c_1\sqrt{T-t}}\int_{L_t}\theta^2\phi\rho|\textbf{K}|^2d\mu_t\nonumber\\
    &&+\frac{c_2 e^{c_1\sqrt{T-t}}}{(T-t)^{\frac{3}{4}}}+c_3 e^{c_1\sqrt{T-t}}.
  \end{eqnarray}
Using the fact that $|\nabla\psi|\leq |\overline{\nabla}\psi|\leq
C$, $|\nabla \theta|=|\textbf{K}|$ and H\"older inequality, we have
\begin{eqnarray}\label{e3.18}
  \int_{L_t}\langle\nabla\psi,\nabla\theta\rangle\phi\rho
   \leq C\left(\int_{L_t}|\textbf{K}|^2\phi\rho\right)^{\frac{1}{2}}\left(\int_{L_t}\phi\rho\right)^{\frac{1}{2}}
   \leq\varepsilon\int_{L_t}|\textbf{K}|^2\phi\rho+C(\varepsilon).
\end{eqnarray}
Combining (\ref{monotheta}) with (\ref{e3.18}) with $\varepsilon$
small enough yields
  \begin{eqnarray}\label{e3.19}
     \frac{\partial}{\partial t}\left(e^{c_1\sqrt{T-t}}\int_{L_t}\theta^2\phi\rho d\mu_t\right)
    &\leq&-e^{c_1\sqrt{T-t}}\int_{L_t}\phi\rho\left(\theta^2
    \left|\textbf{K}+\frac{(F-\textbf{X}_0)^\perp}{2(T-t)}\right|^2+|\textbf{K}|^2\right)d\mu_t\nonumber\\
        &&+\frac{c_2 e^{c_1\sqrt{T-t}}}{(T-t)^{\frac{3}{4}}}+c_3 e^{c_1\sqrt{T-t}}.
  \end{eqnarray}
We will denote $C$ a constant depending on $F_0$, $c_i$ and $T$, which may be different from line to line. First note that, by (\ref{e3.18}),
\begin{equation*}
    \frac{\partial}{\partial t}\left(e^{c_1\sqrt{T-t}}\int_{L_t}\theta^2\phi\rho d\mu_t\right)
     \leq \frac{C}{(T-t)^{\frac{3}{4}}}+C \leq \frac{\partial}{\partial t}\left(-C(T-t)^{\frac{1}{4}}+Ct\right),
\end{equation*}
i.e.,
\begin{equation*}
    \frac{\partial}{\partial t}\left(e^{c_1\sqrt{T-t}}\int_{L_t}\theta^2\phi\rho d\mu_t+ C(T-t)^{\frac{1}{4}}-Ct\right)
     \leq  0.
\end{equation*}
As $e^{c_1\sqrt{T-t}}\int_{L_t}\theta^2\phi\rho d\mu_t+ C(T-t)^{\frac{1}{4}}+Ct$ is bounded for $0<t< T$, we know easily that the limit
\begin{equation*}
    \lim_{t\to T}\left(e^{c_1\sqrt{T-t}}\int_{L_t}\theta^2\phi\rho d\mu_t\right)
\end{equation*}
exists. We denote $\phi_i(F_i(x,s))=\phi(F(x, T+\lambda_i^{-2}s))$. It is easy to see that
\begin{eqnarray*}
 &   & \int_{L^{i}_{s}}\theta_i^2\phi_{i}(F_{i})\frac{1}{0-s}\exp\left(-\frac{|F_{i}|^{2}}{4(0-s)}\right)d\mu ^{i}_{s}\nonumber\\
 & = & \int_{L_{T+\lambda_{i}^{-2}s}}\theta^2\phi(F)\frac{1}{T-(T+\lambda_{i}^{-2}s)}
            \exp\left(-\frac{|F(x,T+\lambda_{i}^{-2}s)-\textbf{X}_0|^{2}}{4(T-(T+\lambda_{i}^{-2}s))}\right)d\mu_{s},
\end{eqnarray*}
where $\phi$ is the function defined in the definition of $\Phi$. Notice that $T+\lambda_{i}^{-2}s\to T$ for any fixed $s<0$. This implies that, for any fixed $s_{1}$ and $s_{2}$ with $-\infty <s_{1}<s_{2}< 0$, we have
\begin{eqnarray}\label{e3.20}
 &   & e^{c_{1}\sqrt{T-(T+\lambda_{i}^{-2}s_{2})}}\int_{L^{i}_{s_{2}}}\theta_i^2\phi_{i}\frac{1}{0-s_{2}}
              \exp\left(-\frac{|F_{i}|^{2}}{4(0-s_{2})}\right)d\mu ^{i}_{s_{2}}\nonumber\\
 &   & -e^{c_{1}\sqrt{T-(T+\lambda_{i}^{-2}s_{1})}}\int_{L^{i}_{s_{1}}}\theta_i^2\phi_{i}\frac{1}{0-s_{1}}
              \exp\left(-\frac{|F_{i}|^{2}}{4(0-s_{1})}\right)d\mu ^{i}_{s_{1}}\nonumber\\
 &   & \to 0  \ \ as \ k\to \infty.
\end{eqnarray}
Integrating (\ref{e3.19}) from $s_{1}$ to $s_{2}$, we obtain
\begin{eqnarray*}
 &   & -e^{c_{1}\sqrt{-\lambda_{i}^{-2}s_{2}}}\int_{L^{i}_{s_{2}}}\theta_i^2\phi_{i}\frac{1}{0-s_{2}}
              \exp\left(-\frac{|F_{i}|^{2}}{4(0-s_{2})}\right)d\mu ^{i}_{s_{2}}\nonumber\\
 &   & +e^{c_{1}\sqrt{-\lambda_{i}^{-2}s_{1}}}\int_{L^{i}_{s_{1}}}\theta_i^2\phi_{i}\frac{1}{0-s_{1}}
              \exp\left(-\frac{|F_{i}|^{2}}{4(0-s_{1})}\right)d\mu ^{i}_{s_{1}}\nonumber\\
 & \geq & \int^{s_{2}}_{s_{1}}e^{c_{1}\sqrt{-\lambda_{i}^{-2}s}}\int_{L^{i}_{s}}\theta_i^2\phi_{i}\rho(F_i,t)
              \left|\textbf{K}_i+\frac{(F_{i})^{\perp}}{2(0-s)}\right|^2d\mu_s^i ds\nonumber\\
 &   & +\int^{s_{2}}_{s_{1}}e^{c_{1}\sqrt{-\lambda_{i}^{-2}s}}\int_{L^{i}_{s}}\phi_{i}\rho(F_i,t)|\textbf{K}_i|^2d\mu_s^i ds\nonumber\\
 &   &    -4c_2\lambda_{i}^{-\frac{1}{2}}((-s_1)^{\frac{1}{4}}-(-s_2)^{\frac{1}{4}})e^{c_{1}\lambda_{i}^{-1}\sqrt{-s_1}}
            -c_{3}\lambda_{i}^{-2}(s_2-s_1)e^{c_{1}\lambda_{i}^{-1}\sqrt{-s_1}}.
\end{eqnarray*}
Thus, we know that
\begin{equation}\label{e3.21}
    \lim_{i\to\infty}\int^{s_{2}}_{s_{1}}e^{c_{1}\sqrt{-\lambda_{i}^{-2}s}}\int_{L^{i}_{s}}\phi_{i}\rho(F_i,t)|\textbf{K}_i|^2d\mu_s^i ds=0.
\end{equation}
In particular, for any $s_1<s_2<0$ and for any $R>0$, we have,
\begin{equation}\label{e3.22}
    \lim_{i\to\infty}\int^{s_{2}}_{s_{1}}\int_{L^{i}_{s}\cap B_R(0)}|\textbf{K}_i|^2d\mu_s^i ds=0.
\end{equation}
This proves the third term of the lemma.

\vspace{.1in}

For the second term, recall that from (\ref{e3.16}), we have
\begin{equation*}
    |\textbf{H}_i|\leq|\textbf{K}_i|+\lambda_i^{-1}n||\overline{\nabla}\psi||_{C^0}\leq |\textbf{K}_i|+C\lambda_i^{-1}.
\end{equation*}
Combining with (\ref{e3.21}), (\ref{e3.22}) and the fact that $\lambda_i\to0$, we obtain
\begin{equation}\label{e3.23}
    \lim_{i\to\infty}\int^{s_{2}}_{s_{1}}e^{c_{1}\sqrt{-\lambda_{i}^{-2}s}}\int_{L^{i}_{s}}\phi_{i}\rho(F_i,t)|\textbf{H}_i|^2d\mu_s^i ds=0,
\end{equation}
and
\begin{equation}\label{e3.24}
    \lim_{i\to\infty}\int^{s_{2}}_{s_{1}}\int_{L^{i}_{s}\cap B_R(0)}|\textbf{H}_i|^2d\mu_s^i ds=0.
\end{equation}

\vspace{.1in}

To prove the first term, we take $f_t\equiv 1$ in (\ref{mono}) to obtain
  \begin{eqnarray}\label{mono1}
     \frac{\partial}{\partial t}\left(e^{c_1\sqrt{T-t}}\int_{L_t}\phi\rho
     d\mu_t\right)
    &\leq&-e^{c_1\sqrt{T-t}}\int_{L_t}\phi\rho
    \left|\textbf{K}+\frac{(F-\textbf{X}_0)^\perp}{2(T-t)}\right|^2d\mu_t\nonumber\\
    &&+\varepsilon
    e^{c_1\sqrt{T-t}}\int_{L_t}\phi\rho|\textbf{K}|^2d\mu_t\nonumber\\
    &&+\frac{c_2 e^{c_1\sqrt{T-t}}}{(T-t)^{\frac{3}{4}}}+c_3 e^{c_1\sqrt{T-t}}.
  \end{eqnarray}
From (\ref{e3.19}) and the above argument, we see that
\begin{equation}\label{Kbdd}
    \int_{0}^{T}e^{c_1\sqrt{T-t}}\int_{L_t}\phi\rho|\textbf{K}|^2d\mu_tdt\leq C<\infty.
\end{equation}
Set $h(t)=\int_{0}^{t}e^{c_1\sqrt{T-u}}\int_{L_u}\phi\rho|\textbf{K}|^2d\mu_udu$ which is a bounded function on $[0,T)$, then arguing as above, we have that
\begin{equation*}
    \frac{\partial}{\partial t}\left(e^{c_1\sqrt{T-t}}\int_{L_t}\phi\rho-\varepsilon h(t) +C(T-t)^{\frac{1}{4}}-Ct\right) d\mu_t
     \leq  0.
\end{equation*}
From this, we conclude that the limit
\begin{equation*}
    \lim_{t\to T}\left(e^{c_1\sqrt{T-t}}\int_{L_t}\phi\rho d\mu_t\right)
\end{equation*}
exists. Arguing in the same way as above to integrating the monotonicity inequality (\ref{e3.24}), we can prove that
\begin{equation*}
    \lim_{i\to\infty}\int^{s_{2}}_{s_{1}}e^{c_{1}\sqrt{-\lambda_{i}^{-2}s}}
    \int_{L^{i}_{s}}\phi_{i}\rho(F_i,t)\left|\textbf{K}_i+\frac{F_i^{\perp}}{2(-s)}\right|^2d\mu_s^i ds=0.
\end{equation*}
Combining with (\ref{e3.21}) yields
\begin{equation}\label{e3.26}
    \lim_{i\to\infty}\int^{s_{2}}_{s_{1}}e^{c_{1}\sqrt{-\lambda_{i}^{-2}s}}
    \int_{L^{i}_{s}}\phi_{i}\rho(F_i,t)\left|\frac{F_i^{\perp}}{2(-s)}\right|^2d\mu_s^i ds=0.
\end{equation}
In particular, for any $s_1<s_2<0$ and for any $R>0$, we have,
\begin{equation}\label{e3.27}
    \lim_{i\to\infty}\int^{s_{2}}_{s_{1}}\int_{L^{i}_{s}\cap B_R(0)}|F_i^{\perp}|^2d\mu_s^i ds=0.
\end{equation}
This finishes the proof the the lemma.
\hfill Q.E.D.

\vspace{.1in}

The following upper bound on volume density is a consequence of monotonicity formula. Similar argument appears in Proposition 2.3 of \cite{CL2}.

\begin{lemma}\label{areaincrease}
  Let $F_t: L\rightarrow M$ be a smooth mean curvature flow of
  a compact Lagrangian submanifold $L_0$ in a compact almost
  Calabi-Yau manifold $M$ of complex dimension $n$. Suppose that
  $L_0$ is zero-Maslov in $M$. For any $\lambda$, $R>0$ and any $s<0$,
  \begin{eqnarray}
    \mu_s^\lambda(L_s^\lambda\cap B_R(0))\leq CR^n,
  \end{eqnarray}
  where $B_R(0)$ is a metric ball in $\textbf{R}^{2n}$ and $C>0$ is
  independent of $\lambda$.
\end{lemma}
{\it Proof.} Set
\begin{eqnarray*}
  \Psi(\textbf{X}_0,t_0,t)=\int_{L_t}\phi(F)\rho(F,\textbf{X}_0
  ,t,t_0)d\mu_t
\end{eqnarray*}
Straightforward computation shows

\begin{eqnarray}\label{mus}
\mu_s^\lambda(L_s^\lambda\cap
B_R(0))&=&\lambda^n\int_{L_{T+\lambda^{-2}s}\cap
B_{\lambda^{-1}}R(X_0)}d\mu_t\nonumber\\
&=&R^n(\lambda^{-1}R)^{-n}\int_{L_{T+\lambda^{-2}s}\cap
B_{\lambda^{-1}R}(X_0)}d\mu_t\nonumber\\
&\leq&CR^n\int_{L_{T+\lambda^{-2}t}\cap
B_{\lambda^{-1}R}(X_0)}\frac{1}{(4\pi)^{n/2}(\lambda^{-1}R)^n}e^{-\frac{|\textbf{X}-\textbf{X}_0|^2}{4(\lambda^{-1}R)^2}}d\mu_t\nonumber\\
&\leq&CR^n\Psi(\textbf{X}_0,T+(\lambda^{-1}R)^2+\lambda^{-2}s,T+\lambda^{-2}s).
\end{eqnarray}
By (\ref{mono1}), we have
\begin{eqnarray}\label{Psi}
  &&\Psi(\textbf{X}_0,T+(\lambda^{-1}R)^2+\lambda^{-2}s,T+\lambda^{-2}s)\nonumber\\
  &\leq&
  e^{c_1\sqrt{T/2+\lambda^{-2}s}}\Psi(\textbf{X}_0,T+(\lambda^{-1}R)^2+\lambda^{-2}s,\frac{T}{2})\nonumber\\
  &&+e^{c_1\sqrt{T/2+\lambda^{-2}s}}\left(\varepsilon\int_{T/2}^{T+\lambda^{-2}s}\int_{L_t}\phi\rho|\textbf{K}|^2d\mu_tds+4c_2(T/2+\lambda^{-2s})^{1/4}+c_3(T/2+\lambda^{-2}s)\right).\nonumber\\
  \end{eqnarray}
Putting (\ref{Psi}) into (\ref{mus}) and using (\ref{Kbdd}),
(\ref{area}), we get
  \begin{eqnarray}
     \mu_s^\lambda(L_s^\lambda\cap
B_R(0))&\leq&CR^n\Psi(\textbf{X}_0,T+(\lambda^{-1}R)^2+\lambda^{-2}s,\frac{T}{2})+CR^n\nonumber\\
&\leq&\frac{\mu_{T/2}(\Sigma_{T/2})}{T^{n/2}}CR^n+CR^n \nonumber\\
&\leq&CR^n.
  \end{eqnarray}
This proves the lemma. \hfill Q.E.D.

\vspace{.1in}

{\it Proof of the Main Theorem.} We follow the argument of the proof of Main Theorem A in \cite{N}.
Pick $s_1<0$ for which
\begin{eqnarray*}
  \lim_{i\rightarrow\infty}\int_{L_{s_1}^i\cap
  B_R(0)}(|F_i^\perp|^2+|\textbf{H}_i|^2+|\textbf{K}_i|^2)=0
\end{eqnarray*}
for all positive $R$.

The maximum principle implies that the Lagrangian angle $\theta_t$
is uniformly bounded and hence, by scale invariance, the same is
true for the Lagrangian angle of $L_s^i$. Lemma \ref{areaincrease}
implies the existence of a constant $D_0$ for which
\begin{eqnarray*}
  \mu_{s_1}^{\lambda_i}(L_{s_1}^i\cap B_R(0))\leq D_0 R^n
\end{eqnarray*}
for all positive $R$. We can, therefore, apply Proposition \ref{Li} to the sequence $(L^i_{s_1})$ and, after a diagonalization argument, obtain a subsequence for which there are integral Special Lagrangian currents
$$L_1, ..., L_N$$
and a finite set $\{\bar{\theta}_1, \cdots, \bar{\theta}_N\}$ such that, for every smooth
function $\phi$ compactly supported, every $f$ in $C^2(\textbf{R})$,
\begin{equation*}
    \lim_{i\to\infty}\int_{L^{i}_{s_1}}f(\theta_{i,s_1})\phi d\mu^{i}_{s_1}=\sum_{j=1}^{N}m_jf(\bar{\theta}_j)\mu_j(\phi),
\end{equation*}
where $\mu_j$ and $m_j$ denote the Radon measure of the support of $L_j$ and its multiplicity respectively.
By our choice of $s_1$, the fact that
\begin{eqnarray*}
  \lim_{i\rightarrow\infty}\int_{L_{s_1}^i\cap B_R(0)}|F_i^\perp|^2d\mu^i_{s_1}=0
\end{eqnarray*}
for all positive $R$ implies that the Special Lagrangians $L_j$ are all cones.

Next, we will show that for all $s_2<0$,
\begin{equation*}
    \lim_{i\to\infty}\int_{L^{i}_{s_2}}f(\theta_{i,s_2})\phi d\mu^{i}_{s_2}
    =\lim_{i\to\infty}\int_{L^{i}_{s_1}}f(\theta_{i,s_1})\phi d\mu^{i}_{s_1}=\sum_{j=1}^{N}m_jf(\bar{\theta}_j)\mu_j(\phi).
\end{equation*}
In fact, from (\ref{e1.3}), (\ref{e2.1}) and the equation (\ref{flow}), we have
\begin{eqnarray*}
% \nonumber to remove numbering (before each equation)
  \frac{d}{ds}\int_{L^{i}_{s}}f(\theta_{i,s})\phi d\mu^{i}_{s}
   &=& \int_{L^{i}_{s}}f'(\theta_{i,s}) \Delta \theta_{i,s}\phi d\mu^{i}_{s}
       +\int_{L^{i}_{s}}f'(\theta_{i,s}) n\lambda_i^{-1}\langle \overline{\nabla}\psi ,\nabla\theta_{i,s}\rangle\phi d\mu^{i}_{s} \\
   & & + \int_{L^{i}_{s}}f(\theta_{i,s})\langle \overline{\nabla}\phi ,\textbf{K}_i\rangle d\mu^{i}_{s}
       -\int_{L^{i}_{s}}f(\theta_{i,s})\phi\langle\textbf{K}_i,\textbf{H}_i\rangle d\mu_s^i.
\end{eqnarray*}
After integration with respect the $s$ variable, all the terms on the tight hand side vanish as $i$ goes to infinity by Lemma \ref{lem3.3}.
We only check the first term. Integrating by parts (and assuming $s_1<s_2$) yields
\begin{eqnarray*}
% \nonumber to remove numbering (before each equation)
  \int_{s_1}^{s_2}\int_{L^{i}_{s}}f'(\theta_{i,s}) \Delta \theta_{i,s}\phi d\mu^{i}_{s}ds
   &=& -\int_{s_1}^{s_2}\int_{L^{i}_{s}}f''(\theta_{i,s}) |\nabla \theta_{i,s}|^2\phi d\mu^{i}_{s}ds \\
   &=& -\int_{s_1}^{s_2}\int_{L^{i}_{s}}f'(\theta_{i,s})\langle \overline{\nabla}\phi ,\nabla\theta_{i,s}\rangle d\mu^{i}_{s}ds
\end{eqnarray*}
Using (\ref{e1.2}) and H\"older's inequality, we have
\begin{equation*}
    \int_{s_1}^{s_2}\int_{L^{i}_{s}}|f''(\theta_{i,s}) |\nabla \theta_{i,s}|^2\phi |d\mu^{i}_{s}ds
    \leq C\int_{s_1}^{s_2}\int_{L^{i}_{s}\cap supp\phi}|\textbf{K}_i|^2 d\mu^{i}_{s}ds
\end{equation*}
and
\begin{equation*}
   \int_{s_1}^{s_2}\int_{L^{i}_{s}}f'(\theta_{i,s})\langle \overline{\nabla}\phi ,\nabla\theta_{i,s}\rangle d\mu^{i}_{s}ds
    \leq C\left(\int_{s_1}^{s_2}\int_{L^{i}_{s}\cap supp\phi}|\textbf{K}_i|^2 d\mu^{i}_{s}ds\right)^{\frac{1}{2}},
\end{equation*}
for some constant $C$ depends on $\phi, f, D_0, s_1, s_2$.

Finally, we show that $\{\bar{\theta}_1, \cdots, \bar{\theta}_N\}$ does not depend on the sequence of rescalings chosen. Let
\begin{equation*}
    (\hat{L}^{k}_{s})_{s<0}
\end{equation*}
be another sequence of rescaled flows for which there are Special Lagrangian cones
$$\hat{L}_1, ..., \hat{L}_N$$
and a finite set $\{\hat{\theta}_1, \cdots, \hat{\theta}_P\}$ such that, for every smooth
function $\phi$ compactly supported, every $f$ in $C^2(\textbf{R})$, and every $s<0$
\begin{equation*}
    \lim_{k\to\infty}\int_{\hat{L}^{k}_{s}}f(\theta_{k,s})\phi d\mu^{k}_{s}=\sum_{j=1}^{P}\hat{m}_jf(\hat{\theta}_j)\hat{\mu}_j(\phi),
\end{equation*}
where $\mu_j$ and $m_j$ denote the Radon measure of the support of $L_j$ and its multiplicity respectively.

For any real number $y$ and any positive integer $q$, we have the
following evolution equation
\begin{eqnarray*}
  \frac{\partial}{\partial
  t}(\theta_t-y)^{2q}=\Delta(\theta_t-y)^{2q}-2q(2q-1)(\theta_t-y)^{2q-2}|\textbf{K}|^2+2qn(\theta_t-y)^{2q-1}\langle\nabla\psi,\nabla
  \theta_t\rangle
\end{eqnarray*}
Take $f_t=(\theta_t-y)^{2q}$ in (\ref{mono}), we get that
\begin{eqnarray*}
  &&\frac{d}{dt}\left(e^{c_1\sqrt{t_0-t}}\int_{L_t}(\theta_t-y)^{2q}\phi\rho
  d\mu_t\right)\\&\leq&
  -e^{c_1\sqrt{t_0-t}}\int_{L_t}2q(2q-1)(\theta_t-y)^{2q-2}\phi\rho|\textbf{K}|^2
  +2qne^{c_1\sqrt{t_0-t}}\int_{L_t}(\theta_t-y)^{2q-1}\langle\nabla\psi,\nabla
  \theta_t\rangle\phi\rho\\&&
  +\varepsilon
  e^{c_1\sqrt{t_0-t}}\int_{L_t}(\theta_t-y)^{2q}\phi\rho|\textbf{K}|^2
  +\frac{c_2e^{c_1\sqrt{t_0-t}}}{(t_0-t)^{\frac{3}{4}}}+c_3e^{c_1\sqrt{t_0-t}}\\.
\end{eqnarray*}

Using H\"{o}lder inequality, the boundedness of $\theta_t$  and the
fact that $|\nabla\theta_t|=|\textbf{K}|$, we have
\begin{eqnarray*}
    &&\frac{d}{dt}\left(e^{c_1\sqrt{t_0-t}}\int_{L_t}(\theta_t-y)^{2q}\phi\rho
  d\mu_t\right)\\&\leq&-e^{c_1\sqrt{t_0-t}}\int_{L_t}2q(2q-1)(\theta_t-y)^{2q-2}\phi\rho|\textbf{K}|^2
  +C\varepsilon
  e^{c_1\sqrt{t_0-t}}\int_{L_t}(\theta_t-y)^{2q-2}\phi\rho|\textbf{K}|^2\\
  &&+C(\varepsilon)e^{c_1\sqrt{t_0-t}}\int_{L_t}(\theta_t-y)^{2q}|\nabla\psi|^2\phi\rho
  +\frac{c_2e^{c_1\sqrt{t_0-t}}}{(t_0-t)^{\frac{3}{4}}}+c_3e^{c_1\sqrt{t_0-t}}\\
\end{eqnarray*}
Choosing $\varepsilon$ small enough, we get
\begin{eqnarray*}
    &&\frac{d}{dt}\left(e^{c_1\sqrt{t_0-t}}\int_{L_t}(\theta_t-y)^{2q}\phi\rho
  d\mu_t\right)\\&\leq&-e^{c_1\sqrt{t_0-t}}\int_{L_t}q(2q-1)(\theta_t-y)^{2q-2}\phi\rho|\textbf{K}|^2\\
  &&+C(\varepsilon)e^{c_1\sqrt{t_0-t}}\int_{L_t}(\theta_t-y)^{2q}|\nabla\psi|^2\phi\rho
  +\frac{c_2e^{c_1\sqrt{t_0-t}}}{(t_0-t)^{\frac{3}{4}}}+c_3e^{c_1\sqrt{t_0-t}}.
\end{eqnarray*}
It follows that (we choose $t_0=T$)
\[\lim_{t\rightarrow T}\left(e^{c_1\sqrt{T-t}}\int_{L_t}(\theta_t-y)^{2q}\phi\rho
  d\mu_t\right)\]
exists, which implies that
\[\lim_{t\rightarrow T}\int_{L_t}(\theta_t-y)^{2q}\phi\rho d\mu_t\]
exists. Thus, by scale invariance, we obtain for any $s,\bar{s}<0$
\begin{eqnarray*}
% \nonumber to remove numbering (before each equation)
  \lim_{i\to\infty}\int_{L_s^i}(\theta_{i,s}-y)^{2q}\phi\rho d\mu_s^i
   &=& \lim_{k\to\infty}\int_{\hat{L}^k_{\bar{s}}}(\theta_{k,\bar{s}}-y)^{2q}\phi\rho d\mu_{\bar{s}}^k \\
   &=& \lim_{t\rightarrow T}\int_{L_t}(\theta_t-y)^{2q}\phi\rho d\mu_t.
\end{eqnarray*}
Therefore,
\begin{equation*}
    \sum_{j=1}^{N}m_j(\bar{\theta}_j-y)^{2q}\mu_j(\phi)=\sum_{j=1}^{P}\hat{m}_j(\hat{\theta}_j-y)^{2q}\hat{\mu}_j(\phi)
\end{equation*}
for any real number $y$ and any positive integer $q$. This implies that
\begin{equation*}
    \{\bar{\theta}_1, \cdots, \bar{\theta}_N\}=\{\hat{\theta}_1, \cdots, \hat{\theta}_P\}.
\end{equation*}
This finishes the proof of the Main Theorem.
\hfill Q.E.D.

\vspace{.1in}

{\it Proof of Corollary \ref{sing}.} Let $\textbf{X}_0$ be a Type I
singularity at $T<\infty$. Choose $t_i\to T$ and $\lambda_i=max_{L\times [0,t_i]}|\textbf{A}|^2$.
Then $\lambda_i$ goes to infinity as $i$ goes to infinity. As $(\textbf{X}_0,T)$ is a Type I
singularity, it is easy to see that the blow up limits $L_\infty$ obtained by the Main Theorem is a smooth minimal Lagrangian
submanifold in $\textbf{C}^n$. Because
$L_\infty$ is smooth, (\ref{limit}) implies $F_\infty^\perp\equiv0$
everywhere. From the monotonicity formula for minimal submanifold of Euclidean space (for example, see Proposition 1.8 of \cite{CM}),
we know that $R^{-n}\mu(L_\infty\cap B_R(0))$ is a constant independent of $R$,
and the volume density ratio at $0$ is one due to the smoothness of
$L_\infty$, so $L_\infty$ is a flat linear subspace of
$\textbf{R}^{2n}$. But the second fundamental form of $L_\infty$ has
length one at 0 according to the blow-up process. This gives the desired contradiction.
\hfill Q.E.D.

\end{document}